\documentclass{amsart}
\usepackage{amscd}

\DeclareMathOperator{\Pic}{Pic} \DeclareMathOperator{\Hom}{Hom}

\def\cA{\mathcal A}

\def\cE{{\mathcal E}}
\def\cF{{\mathcal F}}

\def\cI{{\mathcal I}}

\def\cK{{\mathcal K}}
\def\cL{{\mathcal L}}

\def\cO{{\mathcal O}}

\def\cQ{{\mathcal Q}}

\def\cT{{\mathcal T}}

\def\cV{{\mathcal V}}
\def\cW{{\mathcal W}}

\def\M{{\mathfrak M}}

 \newtheorem{thm}{Theorem}[section]
 \newtheorem{cor}[thm]{Corollary}
 \newtheorem{lem}[thm]{Lemma}
 \newtheorem{prop}[thm]{Proposition}
 \newtheorem{defn}[thm]{Definition}

 \newcommand{\map}{\longrightarrow }

\begin{document}

\title[The number of rational curves on K3 surfaces]
 {The number of rational curves on K3 surfaces}

\author{Baosen Wu}

\address{Department of Mathematics, Stanford University, Stanford, CA 94305}

\email{bwu@math.stanford.edu}


%%% ----------------------------------------------------------------------
\begin{abstract}

Let $X$ be a $K3$ surface with a primitive ample divisor $H$, and
let $\beta=2[H]\in H_2(X, \mathbf Z)$. We calculate the
Gromov-Witten type invariants $n_{\beta}$ by virtue of Euler
numbers of some moduli spaces of stable sheaves. Eventually, it
verifies Yau-Zaslow formula in the non primitive class $\beta$.

\end{abstract}

%%% ----------------------------------------------------------------------
\maketitle

%%% ----------------------------------------------------------------------
\section* {Introduction}

Let $X$ be a $K3$ surface with an ample divisor $H$, and let $C\in |H|$ be a reduced curve. By adjunction
formula, the arithmetic genus of $C$ is $g=\frac{1}{2}H^2+1$. Under the assumption that the homology class
$[H]\in H_2(X, \mathbf Z)$ is primitive, Yau and Zaslow \cite{YZ} showed that the number of rational curves in
the linear system $|H|$ is equal to the coefficient of $q^g$ in the series
\begin{align}
\frac {q}{ \Delta (q)}=\prod_{k>0}\frac{1}{(1-q^k)^{24}} \notag &= \sum_{d\ge 0} G_dq^d\\ \notag &=
1+24q+324q^2+3200q^3+25650q^4+176256q^5+\cdots \notag
\end{align}
Here a multiplicity $e(\bar JC)$ is assigned to each rational
curve $C$ in the counting(\cite{B}).

In \cite{FGS}, Fantechi, G\"ottsche and van Straten gave an interpretation of the multiplicity $e(\bar JC)$. Let
$M_{0, 0}(X,[H])$ be the moduli space of genus zero stable maps $f:\mathbf P^1\to X$ with $f_*([\mathbf
P^1])=[H]\in H_2(X, \mathbf Z)$. $M_{0, 0}(X,[H])$ is a zero dimensional scheme which is in general nonreduced.
Let $\iota:C\hookrightarrow X$ be a rational curve in the class $[H]$, and $n:\mathbf P^1\to C$ its
normalization. Then $f=\iota\circ n:\mathbf P^1\to X$ is a closed point of $M_{0, 0}(X,[H])$ and $e(\bar JC)$ is
equal to the multiplicity of $M_{0, 0}(X,[H])$ at $f$.

There is another formulation and generalization of Yau and
Zaslow's formula by virtue of Gromov-Witten invariants. For $K3$
surfaces, the usual genus $0$ Gromov-Witten invariants vanish. To
remedy this, one can use the notion of twistor family developed by
Bryan and Leung in \cite{BL} provided that $\beta$ is a primitive
class. In general, there is an algebraic geometric approach
proposed by Jun Li \cite{L} using virtual moduli cycles. Roughly
speaking, he defines Gromov-Witten type invariants $N_g(\beta)$ on
$K3$ surfaces by modifying the usual tangent-obstruction complex.
When $\beta$ is primitive, these invariants coincide with those
defined by twistor family. Geometrically, $N_g(\beta)$ can be
thought as Gromov-Witten invariants of a one dimensional family of
$K3$ surfaces, which actually count curves in the original
surface. For the rigorous definitions, see \cite{BL},\cite{L}.

Bryan and Leung \cite{BL} proved a formula for $N_g(\beta)$ when
$\beta$ is primitive. Let $n_{\beta}=N_0(\beta)$. Then
$n_{\beta}=G_d$ with $d=\frac{1}{2}\beta^2+1$. It recovers the
formula of Yau and Zaslow. For a non primitive class $\beta$, the
numbers $N_g(\beta)$ are still unknown. However, there is a
conjectural formula for $N_0(\beta)$(\cite{L}). Using the notation
$n_{\beta}$, it says
\[n_{\beta}=\sum_k\frac{1}{k^3}G_{\frac{1}{2}(\frac{\beta}{k})^2+1}\]
where the sum runs over all integers $k>0$ such that
$\frac{\beta}{k}$ is an integral homology class(see also
\cite{G}). The case $\beta=2[H]$ with $[H]$ primitive and $H^2=2$
was proved by Gathmann in \cite{G}.

In this paper, we will prove the following result.
\begin{thm}
{Let $X$ be a K3 surface with an ample divisor $H$. Assume $[H]\in H_2(X,\mathbf Z)$ is primitive. Let
$\beta=2[H]$ and $g=\frac{1}{2}H^2+1$. Then
$$n_{\beta}=G_{4g-3}+\frac{1}{8}G_g.$$
}\end{thm}

Now we sketch the proof of this theorem. It can be divided into two steps. First, we deform the pair $(X, H)$ to general position and then reduce the
calculation of $n_{\beta}$ to $N_{\beta}$, which is the number of reduced and irreducible rational curves in $\beta$. The second step is the
calculation of $N_{\beta}$. In Gathmann's approach \cite{G}, the assumption $H^2=2$ is essential in this step. In this paper, we will generalize the
approach of Yau and Zaslow \cite{YZ} according to the suggestion in \cite{L}.

Next, we describe these two steps in details.

We begin with the first step. Let $(X, H)$ be a pair of a $K3$
surface $X$ and a primitively polarization $H$ on $X$. It is well
known that two pairs $(X, H)$ and $(X', H')$ with $H^2=H'^2$ are
deformation equivalent. One can choose a general primitively
polarized $K3$ surface $(X, H)$, such that $\Pic X=\mathbf Z\cdot
[H]$ and every rational curve in the linear system $|H|$ is nodal
\cite{C}. Moreover, using a generalization of the method in
\cite{C}, one can also assume that any two rational curves in the
system $|H|$ intersect transversely \cite{C2}. Now we fix such a
pair $(X, H)$ once and for all. Since $n_{\beta}$ is a deformation
invariant, we only need to calculate $n_{\beta}$ for such a
surface.

By the enumerative interpretation of $n_{\beta}$, and follow up a similar argument as in \cite{G}, all stable
maps $f: C\to X$ with $f_*([C])=\beta$ can be decomposed into the following three types:

1) The domain $C$ is $\mathbf P^1$, and the image $f(C)\subset X$ is a reduced and irreducible curve in the
linear system $|2H|$. We denote the number of such maps by $N_{\beta}$. The multiplicity of such $f$ is the
Euler number of the compactified Jacobian of the image $f(C)$, as shown in \cite{B}.

2) The domain is a union of two $\mathbf P^1$ that intersect at
one point $P$. In this case, the image is a union of two rational
nodal curves that intersect at $H^2$ points. The image of $P$ has
to be one of the intersections, hence there are $H^2$ such maps.
Since the number of rational curves in the system $|H|$ is $G_g$,
the total number of such maps is $\frac{1}{2} G_g(G_g-1)H^2$.

3) $f: C\to X$ is a double cover onto the image $f(C)$. There are two different cases:

(a) Double covers that factor through the normalization of $f(C)$, this space has dimension $2$.

(b) Double covers that do not factor through the normalization. In
this case, the domain must be a union of two $\mathbf P^1$, which
intersect at one point $P$. The image of $P$ is a node on the
image curve $f(C)$, and there is only one map for each choice of
node. Note that the number of nodes on $f(C)$ is equal to the
arithmetic genus $g$, so there are totally $gG_g$ such maps.

By Lemma $4.1$ in \cite{G}, the contribution of type (3a) is $\frac{1}{8} G_g$. Therefore,
\begin{align}
n_{\beta}=N_{\beta}+\frac{1}{2} G_g(G_g-1)H^2+gG_g+\frac{1}{8} G_g.
\end{align}

Since the first step of the proof is already known, in this paper,
we will focus on the second step, namely, the calculation of the
number $N_{\beta}$ of reduced and irreducible rational curves in
the linear system $|2H|$. To this end, we will work with the
moduli space of sheaves on a $K3$ surface.

Let $(X,H)$ be the pair we fixed previously. Let $\M$ be the
moduli scheme of stable sheaves $\cF$ on $X$ such that
$\dim\cF=1$, $c_1(\cF)=2H$ and $\chi(\cF)=1$. The Hilbert
polynomial of $\cF$ with respect to the polarization $H$ is
$2H^2\cdot n+1$. Since there is no strictly semistable sheaf in
$\M$, by \cite{M}, $\M$ is a smooth projective variety, and its
Euler number $e(\M)$ is $G_{2H^2+1}$(\cite{Yo}). In section $1$,
we will construct a morphism $\Phi: \M\to |2H|$ that sends
$\cF\in\M$ to its support in $|2H|$. For $D\in |2H|$, we denote by
$\M_D$ the fiber of $\Phi$ over $D$ with the reduced subscheme
structure. When $D$ is reduced and irreducible, $\M_D$ is the
compactified Jacobian $\bar JD$ of $D$. In section $2$, we will
show that $e(\M_D)=0$ if $D$ has an irreducible component whose
geometric genus is positive.

Therefore only divisors with rational components contribute to the
Euler number $e(\M)$. Since $H$ is primitive, we have three types
of these divisors in the linear system $|2H|$.

1) $D=C$, $C$ is a rational curve in homology class
$\beta$($=2[H]$). In this case, $\M_D\cong\bar JD$. The number of
such divisors $D$, counted with multiplicity $e(\bar JD)$, is
equal to $N_{\beta}$.

2) $D=C_1+C_2$, where $C_1$ and $C_2$ are different rational nodal
curves. In this case, both $C_i$ are contained in the linear
system $|H|$. There are totally $\frac{1}{2}G_g(G_g-1)$ divisors
of this type. We will show that $e(\M_D)=H^2$ in section $3$.

3) $D=2C_0$, where $C_0$ is a rational nodal curve and contained
in $|H|$. The number of such divisors is $G_g$. In the last two
sections we will prove $e(\M_D)=g$, which is equal to the number
of nodes of $C_0$.

Since $e(\M)=\sum e(\M_D)$, where the sum runs over all divisors
$D$ with rational components, we get
\begin{align}
N_{\beta} &= e(\M)-\frac{1}{2}G_g(G_g-1)H^2-gG_g \\\notag &=
G_{4g-3}-\frac{1}{2}G_g(G_g-1)H^2-gG_g \notag
\end{align}

Together with (1), we prove
\[n_{\beta}=G_{4g-3}+\frac{1}{8}G_g.\]

Recently, J. Li and the author \cite{LW} proved the conjectured formula for non primitive class $\beta=n[H]$ with $n<6$, under the assumption that
the transversality of rational curves still holds.

I am most grateful to Jun Li, from whom I learned moduli spaces of sheaves and Gromov-Witten invariants. During the preparation of this paper, his
constant encouragement and discussions are invaluable. After finishing the manuscript, the author is informed that J.Lee and N.C.Leung \cite{LL1}
proved the same result using degeneration method and also counted genus $1$ curves in $K3$ surfaces \cite{LL2}.

\section {Decomposition of the moduli scheme $\M$}

We start with some definitions and notations(\cite{Si},\cite{HL}).

Let $X$ be a complex projective scheme with an ample line bundle
$\cO(1)$. For a coherent sheaf $\cE$ of $\cO_X$-module, the
Hilbert polynomial $p(\cE, n)$ of $\cE$ is defined as
\[p(\cE, n)=\dim H^0(X, \cE(n)),\ n\gg 0.\]
The dimension of the support of $\cE$ is equal to the degree of
$p(\cE, n)$. A coherent sheaf $\cE$ is pure of dimension $d$ if
for any nonzero coherent subsheaf $\cF\subset\cE$, $\dim\cF=d$.

The Hilbert polynomial $p(\cE, n)$ can be written as
\[p(\cE, n)=\frac{a_0}{d!}n^d+\frac{a_1}{(d-1)!}n^{d-1}+\cdots \]
with integral coefficients $a_i=a_i(\cE)$. We define the slope of
$\cE$ to be
\[\mu(\cE)=a_0(\cO_X)\frac{a_1(\cE)}{a_0(\cE)}-a_1(\cO_X).\]

\begin{defn}
{A coherent sheaf $\cE$ is stable (resp. semistable) if it is
pure, and if for any nonzero proper subsheaf $\cF\subset\cE$,
there exists an $N$, such that for $n>N$,
$$\frac{p(\cF, n)}{a_0(\cF)}<\frac{p(\cE, n)}{a_0(\cE)}\quad (resp. \le).$$
}\end{defn}

\begin{defn}
{A coherent sheaf $\cE$ is $\mu$-stable (resp. $\mu$-semistable)
if it is pure, and if for any nonzero proper subsheaf $\cF\subset
\cE$,
$$\mu(\cF)<\mu(\cE)\quad (resp. \le).$$
}\end{defn}

\begin{thm}
{\cite{Si}Let $X$ be a complex projective scheme with an ample
line bundle $\cO(1)$. There is a projective coarse moduli scheme
whose closed points represent the $S$-equivalence classes of
semistable sheaves with Hilbert polynomial $P(n)$.  }\end{thm}

Let $X$ be a K3 surface with an ample line bundle $H$. By
Riemann-Roch theorem, the Hilbert polynomial of a torsion free
sheaf $\cE$ is
\[p(\cE, n)=\frac{r}{2}H^2n^2+(c_1\cdot
H)n+r\chi(\cO_X)+\frac{1}{2}(c_1^2-2c_2),\] where $r$ is the rank
of $\cE$ and $c_i=c_i(\cE)$. Let $\cF$ be a pure sheaf of
dimension $1$ on $X$. By a locally free resolution, one can verify
that the Hilbert polynomial of $\cF$ is
\[p(\cF,
n)=(c_1(\cF)\cdot H)n+\frac{1}{2}(c_1^2(\cF)-2c_2(\cF)).\] It is
clear that for such sheaves the notion of stability and
$\mu$-stability coincide.

From now on, we fix a pair $(X, H)$ of a $K3$ surface $X$ and a polarization $H$ of $X$, such that

1) $\Pic X=\mathbf Z\cdot [H]$;

2) every rational curve in $|H|$ is nodal; and

3) any two distinct rational curves in $|H|$ intersect transversely.

We let $\beta=2[H]\in H_2(X,\mathbf Z)$. Our immediate goal is to
calculate $N_{\beta}$, the number of reduced and irreducible
rational curves in $|2H|$ counted with multiplicity. To this end,
we consider the moduli scheme $\M$ of stable sheaves $\cF$ of
$\cO_X$-modules that satisfy $\dim\cF=1$, $c_1(\cF)=\beta$ and
$\chi(\cF)=1$.

\begin{thm}
{\cite{Yo} $\M$ is a smooth projective variety. The Euler number
$e(\M)$ is $G_{2H^2+1}$. }\end{thm}

Next we define the morphism $\Phi:\M\to |2H|$ mentioned earlier.

Let $\cF$ be a sheaf in $\M$. Since $\cF$ is pure of dimension
$1$, it admits a length $1$ locally free resolution
\[ 0\map \cE_1\stackrel{f}{\map}\cE_0\map\cF\map 0, \]
with $r(\cE_1)=r(\cE_0)$. The homomorphism $f:\cE_1\to\cE_0$
induces a homomorphism $\wedge ^rf: \wedge ^r \cE_1\to\wedge ^r
\cE_0$, and a nonzero global section $s\in H^0(\wedge ^r
\cE_1)^{-1}\otimes (\wedge ^r \cE_0))$, that defines an effective
divisor $D=s^{-1}(0)$ on $X$. Since $(\wedge ^r \cE_1)^{-1}\otimes
(\wedge ^r \cE_0)=c_1(\cF)=2H$, $D$ is contained in the linear
system $|2H|$. The assignment $\cF\to D$ defines a morphism $\Phi:
\M\to |2H|$.

We now give a specific decomposition of the projective space
$|2H|$ according to the topological type of $D\in |2H|$.

We let $\cW_1$ be the set of divisors $D$ which is reduced and irreducible. The arithmetic genus of $D$ is $p_a(D)=2H^2+1$, which is an invariant for
all $D\in\cW_1$. We further stratify $\cW_1$ according to the geometric genus of curves, $\cW_1=\sqcup_k \cW_1^k$, where $\cW_1^k$ consists of those
$D$ that have geometric genus $k$. Clearly $\sqcup_{k\le a} \cW_1^k$ is closed in $\cW_1$. Let $\cW_2$ be the stratum of divisors $D=C_1+C_2$ with
$C_1 \ne C_2$ and $C_i\in |H|$. Without loss of generality, we can assume $p_g(C_1)\le p_g(C_2)$. For $a\le b$, we let $\cW_2^{a, b}\subset\cW_2$ be
the subset of divisors $D$ with $p_g(C_1)=a$ and $p_g(C_2)=b$. Then $\cW_2=\sqcup \cW_2^{a, b}$. Let $\cW_3$ be the subset of divisors $D=2C_0$ with
$C_0\in |H|$. Similarly, $\cW_3=\sqcup_k\cW_3^k$, where $\cW_3^k$ consists of $D=2C_0$ with $p_g(C_0)=k$.

Put together, \[|2H|=(\sqcup_k\cW_1^k)\bigsqcup(\sqcup_{a\le
b}\cW_2^{a, b})\bigsqcup(\sqcup_k\cW_3^k).\] This induces a
decomposition on $\M$,
\[\M=(\sqcup_k\Phi^{-1}(\cW_1^k))\bigsqcup(\sqcup\Phi^{-1}(\cW_2^{a,
b}))\bigsqcup(\sqcup_k\Phi^{-1}(\cW_3^k)).\]

Now we state a general fact on the Euler number of varieties.

Let $Z$ be a complex variety. Let $Z=\sqcup Z_i$ be a
decomposition into locally closed subset $Z_i$. Then the Euler
number $e(Z)=\sum e(Z_i)$.

Apply this to the decomposition of $\M$, we have
\begin{align}
e(\M) \notag &=
e(\Phi^{-1}(\cW_1))+e(\Phi^{-1}(\cW_2))+e(\Phi^{-1}(\cW_3))\\
\notag &= \sum_k e(\Phi^{-1}(\cW_1^k))+\sum_{a\le b}
e(\Phi^{-1}(\cW_2^{a, b}))+\sum_k e(\Phi^{-1}(\cW_3^k)). \notag
\end{align}

\begin{prop}
{\cite{B} Let $h: Y\to Z$ be a surjective morphism between complex
algebraic varieties. Suppose that $e(h^{-1}(z))=0$ for every
closed point $z\in Z$. Then $e(Y)=0$. }\end{prop}

The following proposition will be proved in the next section.

\begin{prop}
{Suppose $D$ is a divisor that has one irreducible component whose
geometric genus is positive, then $e(\M_D)=0$. }\end{prop}

Combine these results, we have
\[e(\M)=e(\Phi^{-1}(\cW_1^0))+e(\Phi^{-1}(\cW_2^{0,
0}))+e(\Phi^{-1}(\cW_3^0)).\]

Because $e(\Phi^{-1}(\cW_1^0))$ is equal to $N_{\beta}$, and
$e(\M)=G_{2H^2+1}$, To calculate $N_{\beta}$, it suffices to find
the Euler numbers $e(\Phi^{-1}(\cW_2^{0, 0}))$ and
$e(\Phi^{-1}(\cW_3^0))$. The number $e(\Phi^{-1}(\cW_2^{0, 0}))$
is essentially known, which is equal to $\frac{1}{2}G_g(G_g-1)H^2$
as will be shown in section $3$. The main body of the remainder of
the paper is to show that $e(\Phi^{-1}(\cW_3^0))=gG_g$. Therefore,
\[ N_{\beta}=G_{2H^2+1}-\frac{1}{2}G_g(G_g-1)H^2-gG_g.\]
Apply equality (1) in the introduction, we obtain the formula in
the main theorem.

\section {Proof of Proposition 1.6}

We state a basic fact about the Euler number of a variety. Let $X$
be a quasi-projective variety. If there exists a finite group
action on $X$ which is free of fixed point, then $e(X)$ is
divisible by the order of this group. Therefore, if for any
positive integer $N$, there is a finite group $G_N$ whose order is
greater than $N$, and a free $G_N$ action on $X$, then $e(X)$ is
zero.

If $D$ is a reduced and irreducible curve, then $\M_D\cong\bar JD$. Since the geometric genus of $D$ is positive, $e(\bar JD)=0$(see\cite{B}). Now if
$D=C_1+C_2$ with $C_i\in |H|$ and by assumption the geometric genus $p_g(C_2)>0$. From the restriction homomorphism $\alpha:\Pic D\to \Pic C_2$, we
can choose a subgroup $G\subset \Pic D$, such that for $\cL\in G$, $\cL|_{C_1}\cong \cO_{C_1}$ and $\alpha(\cL)=\cL|_{C_2}$ is trivial if and only if
$\cL$ is trivial. Next we show that the $G$-action on $\frak M_D$ defined by tensorization is free, i.e., for any sheaf $\cF\in\frak M_D$ and $\cL\in
G$, $\cF\otimes \cL\cong\cF$ if and only if $\cL$ is trivial. To this end, suppose $\cF\otimes \cL\cong\cF$ for some $\cF$ and $\cL$. Let $\cF_2$ be
the torsion free part of the restriction $\cF|_{C_2}$. We obtain $\cF_2\otimes\alpha(\cL)\cong \cF_2$ and therefore $\alpha(\cL)$ is trivial by the
same argument as in case 1. Finally, it implies $\cL$ is trivial by the choice of the subgroup $G$. Finally, $D=2C_0$ is a divisor whose associated
subscheme is a nonreduced curve $C$, and a closed point in $\frak M_D$ is a sheaf of $\cO_C$-modules. To prove this case, we first recall some facts
on nonreduced curves.

Let $\cF$ be a sheaf of $\cO_X$-modules. An infinitesimal
extension(\cite{Ha}, Exer II 8.7) of $X$ by $\cF$ is a scheme
$X'$, with an ideal sheaf $\cI$, such that $\cI^2=0$ and $(X',
\cO_{X'}/\cI)\cong (X, \cO_X)$ and such that $\cI$ with the
induced structure of $\cO_X$-module is isomorphic to the given
sheaf $\cF$. Let $S$ be a smooth projective surface, and
$C_0\subset S$ be a reduced and irreducible curve. There is an
associated closed subscheme $C\subset S$ to the divisor $2C_0$. In
fact, $C$ is an infinitesimal extension of $C_0$ by
$\cI=\cO_S(-C_0)|_{C_0}$.

Next we discuss the Picard group of $C$(\cite{Ha}, Exer III 4.6).
From the exact sequence of sheaves of abelian groups
\[0\map\cI\map\cO_C^*\map\cO_{C_0}^*\map 0,\] there is an induced
exact sequence
\[ 0\map H^1(C, \cI)\map \Pic C\map \Pic C_0\map
0.\] Notice that $H^1(C, \cI)$ is a vector space and hence an
injective $\mathbf Z$-module, it implies that $\Pic C\cong \Pic
C_0\oplus H^1(C, \cI)$ as groups. For $\cL\in \Pic C$, we let
$\cL_0\in\Pic C_0$ be the restriction of $\cL$ to $C_0$.

Now we continue the proof. Let $\pi:\tilde C_0\to C_0$ be the
normalization of $C_0$. Then $\Pic^0 C_0\cong
\Pic^0\tilde{C_0}\oplus A$, where $A$ is an affine commutative
group. Since the genus of $\tilde C_0$ is positive, $\Pic ^0\tilde
C_0$ is nontrivial. For any odd prime $p$, we can choose an order
$p$ subgroup $G\subset\Pic^0 C$, such that for $\cL\in G$, $\cL$
is trivial if and only if $\tilde{\cL}=\pi^*\cL_0$ is trivial.
There is a $G$-action on $\M_D$ defined by tensorization. Next we
show that this group action is free and therefore $e(\M_D)=0$.

Suppose $\cL\otimes\cE\cong \cE$ for some sheaf $\cE\in\M_D$ and
$\cL\in G$. Restrict to $C_0$ and let $\cE_0$ be the torsion free
part of $\cE\otimes\cO_{C_0}$, we get $\cL_0\otimes\cE_0\cong
\cE_0$.

1) If $\cI\cE\ne 0$, $\cE_0$ is a rank $1$ torsion free sheaf on
$C_0$. Using the same argument as in case $1$, we obtain
$\cL_0\cong \cO_{C_0}$, and hence $\cL\cong \cO_C$, i.e., the
group action is free.

2) If $\cI\cE=0$, $\cE_0$ is a rank $2$ torsion free sheaf on
$C_0$. Let $\tilde \cE_0$ be the torsion free part of
$\pi^*\cE_0$. Then we have $\tilde\cL_0\otimes\tilde
\cE_0\cong\tilde \cE_0$. Take top wedge on both sides, we get
$\tilde\cL_0^{\otimes
2}\otimes\wedge^2\tilde\cE_0\cong\wedge^2\tilde\cE_0$. Since
$\wedge^2\tilde\cE_0$ is invertible, $\tilde\cL_0^{\otimes
2}\cong\cO_{\tilde{C_0}}$. Note that $\cL\in G$ and $G$ has odd
prime order $p$, it implies that $\cL\cong\cO_C$. Hence the group
action is free.

\section {Calculation of $e(\Phi^{-1}(\cW_2^{0, 0}))$}

Recall that $\cW_2^{0,0}$ is a finite set of divisors $D=C_1+C_2$
with $C_1, C_2\in |H|$ being rational nodal curves and intersect
transversally, $\Phi^{-1}(\cW_2^{0, 0})=\sqcup \M_{D_i}$ with
$D_i\in\cW_2^{0,0}$. We will calculate $e(\M_{D})$ for
$D\in\cW_2^{0,0}$ and then the Euler number $e(\Phi^{-1}(\cW_2^{0,
0}))$ follows.

A closed point in $\M_D$ is a stable sheaf $\cE$ of
$\cO_D$-modules, such that the restrictions $\cE|_{C_i}$ are rank
$1$ sheaves of $\cO_{C_i}$-modules respectively. Let $x_1, x_2,
\cdots, x_s$ be a list of intersections of $C_1$ and $C_2$. Then
$s=H^2>0$. Since $\cE$ is stable, there is at least one point
$x_k$, so that the stalk $\cE_{x_k}$ is isomorphic to $\cO_{x_k}$.
For otherwise, $\cE$ is the direct image of some sheaf on the
disjoint union of $C_1$ and $C_2$, which violates the stability of
$\cE$.

We let $S_{ij}\subset\M_D$ be the subset of stable sheaves $\cE$
such that $\cE_{x_i}\cong\cO_{x_i}$ and $\cE_{x_j}\cong\cO_{x_j}$
for two intersection points $x_i$ and $x_j$. We can find a
subgroup $G\subset\Pic D$ coming from the gluing of $\cO_{C_1}$
and $\cO_{C_2}$ at $x_i$ and $x_j$. $G\cong \mathbf C^*$. Now
follow a similar argument as in the previous section, the
$G$-action on $S_{ij}$ is free. Therefore, the contribution to the
Euler number $e(\M_D)$ come from stable sheaves $\cE$ whose stalks
are not $\cO$ at all nodes but one intersection point. Since both
$C_i$ are rational curves, there is only one such stable sheaf
corresponds to an intersection point. We have
\begin{prop}
{Let $D$ be a divisor in the set $\cW_2^{0,0}$. Then
$e(\M_D)=H^2$. }\end{prop}

Since the number of rational curves in $|H|$ is $G_g$,
$\cW_2^{0,0}$ is a finite set with cardinality
$\frac{1}{2}G_g(G_g-1)$.

\begin{cor}
{$e(\Phi^{-1}(\cW_2^{0, 0}))=\frac{1}{2}G_g(G_g-1)H^2$. }\end{cor}

\section {Calculation of $e(\Phi^{-1}(\cW_3^0))$, Part I}

In the remainder of this paper, we will calculate the Euler number
$e(\Phi^{-1}(\cW_3^0))$. Remember that $\cW_3^0$ is a finite set
of divisors $D=2C_0$ with $C_0\in |H|$ being rational nodal
curves, there is a decomposition $\Phi^{-1}(\cW_3^0)=\sqcup \frak
M_{D_i}$. It suffices to calculate $e(\frak M_D)$ for
$D\in\cW_3^0$.

Recall that for every effective divisor, there is an associated
subscheme. Let $C$ be the nonreduced curve associated to $D=2C_0$.
Then every closed point in $\frak M_D$ corresponds to a stable
sheaf $\cE$ of $\cO_C$-modules, such that the Hilbert polynomial
$P_{\cE}(n)$ is $2H^2+1$.

There are two kinds of these sheaves. A sheaf $\cE$ in the first
type satisfies $\cI\cE=0$, where $\cI\subset\cO_C$ is the
nilpotent ideal sheaf. That is to say, $\cE$ is a rank $2$ sheaf
on $C_0$, the reduced part of $C$. Let $\M_D^1\subset\M_D$ be the
subset of sheaves of this type.  The second type consists of
sheaves $\cE$ satisfy $\cI\cE\ne 0$. It is direct to verify that
for sheaves of this type, $\cE_\eta\cong\cO_\eta$ with $\eta$ the
generic point of $C$. Let $\M_D^2$ be the subset of sheaves of the
second type. Then $e(\M_D)=e(\M_D^1)+e(\M_D^2)$.

In this section, we calculate $e(\frak M_D^1)$. The discussion of $\frak M_D^2$ is left to the next section. The result has been obtained by T.
Teodorescu in his PhD thesis \cite{T} which deals with a more general problem. It is also proved in \cite{W} independently. Now we use a slightly
different approach.

We recall some standard facts about sheaves on a nodal curve.
Since we will not talk about nonreduced curves in the next part of
this section, we use $C$, instead of $C_0$, to denote a nodal
curve. We always work on the complex topology.

Let $C$ be a projective curve with $n$ ordinary nodes $x_1, x_2,
\cdots, x_n$ as singularities, and let $\pi: \tilde C\to C$ be the
normalization of $C$. A torsion free sheaf $\cE$ is locally free
away from the nodes. It has the following nice local structure at
each node $x_i\in C$(\cite{Se})
\[\cE_{x_i}\cong\cO_{x_i}^{\oplus a_i}\oplus m_{x_i}^{\oplus
(r-a_i)},\] where $m_{x_i}\subset\cO_{x_i}$ is the maximal ideal,
and $r$ is the rank of $\cE$. Let $\hat\pi: \hat C \to C$ be a
partial normalization of $C$ at one node $x$. Then there exists a
torsion free sheaf $\cF$ on $\hat C$ such that
$\cE\cong\hat\pi_*\cF$ if and only if $\cE_{x}\cong m_{x}^{\oplus
r}$.

Let $r\ge 1$ be an integer and choose $n$ such that $(r, n)=1$.
There is a smooth projective variety $\M_C(r, n)$ whose closed
points correspond to isomorphism classes of stable $\cO_C$-modules
$\cE$, such that $r(\cE)=r$ and $\chi(\cE)=n$.

Next we introduce the notion of admissible quotients. It will be
used to determine whether two torsion free sheaves $\cE_1$,
$\cE_2$ are isomorphic. Let $\cE$ be a torsion free sheaf and
$(\pi^*\cE)^\sharp$ be the torsion free part of $\pi^*\cE$. There
is a canonical exact sequence
\[ 0\map\cE\map\pi_*(\pi^*\cE)^{\sharp}\map\cT\map 0,\]
where $\cT\cong\oplus\mathbf C_{x_i}^{\oplus a_i}$ is a skyscraper
sheaf supported at the nodes.
\begin{defn}
{Let $\cV$ be a rank $r$ locally free sheaf on $\tilde C$, and
$\cQ=\oplus\mathbf C_{x_i}^{\oplus a_i}$ be a skyscraper sheaf
supported at the set of nodes on $C$. Let $\rho:\pi_*\cV\to\cQ$ be
a surjective morphism and $\cE$ be the kernel of $\rho$. $\rho$ is
said to be an admissible quotient if there is a commutative
diagram
\begin{equation*}
\begin{CD}
0@>>>\cE@>>>\pi_*\cV@>\rho >>\cQ @>>> 0\\
@.       @VV=V        @VV\cong V   @VV\cong V \\
0@>>>\cE@>>>\pi_*(\pi^*\cE)^{\sharp}@>>>\cT@>>> 0\\
\end{CD}
\end{equation*}
where the second row is the canonical exact sequence. }\end{defn}

Let $p_i, q_i\in\tilde C$ be the inverse images of the node $x_i$.
Then $(\pi_*\cV)_{x_i}=\cV_{p_i}\oplus\cV_{q_i}$. The homomorphism
$\rho$ is given by
$\rho_{x_i}:\cV_{p_i}\oplus\cV_{q_i}\to\cQ_{x_i}$. Let
$\iota_i^1:\cV_{p_i}\to \cV_{p_i}\oplus\cV_{q_i}$ and
$\iota_i^2:\cV_{q_i}\to \cV_{p_i}\oplus\cV_{q_i}$ be the natural
injections. By the definition of admissible quotients,
$\rho_{x_i}^k=\rho_{x_i}\circ\iota_i^k$ are both surjective.
Conversely, we have

\begin{prop}
{Let $\rho:\pi_*\cV\to\cQ$ be a quotient such that $\rho_{x_i}^k$
defined above are surjective for all $i=1, 2 \cdots, n$ and $k=1,
2$. Then $\rho$ is an admissible quotient. }\end{prop}

\begin{proof}
{Let $\cE$ be the kernel of $\rho$. Apply the functor $\pi^*$ to
the exact sequence
\[ 0\map\cE\map\pi_*\cV\stackrel{\rho}{\map}\cQ\map 0, \]
we get
\[\pi^*\cE\stackrel{\psi_1}{\map}\pi^*(\pi_*\cV)\stackrel{\psi_2}{\map}\pi^*\cQ\map 0.\]
Since $\rho^k_{x_i}$ are surjective, the restriction of $\psi_2$
on the torsion part $\cT'\subset\pi^*(\pi_*\cV)$ is surjective,
i.e. $\psi_2(\cT')=\pi^*\cQ$. It implies that the homomorphism
$\pi^*\cE\to(\pi^*(\pi_*\cV))^{\sharp}$induced by $\psi_1$ is
surjective. Because the kernel of $\psi_1$ is a skyscraper sheaf,
$\psi_1$ induces an isomorphism
$(\pi^*\cE)^{\sharp}\to(\pi^*(\pi_*\cV))^{\sharp}$. Since every
step is functorial, the result follows from the canonical
isomorphism $(\pi^*(\pi_*\cV))^{\sharp}\cong\cV$. }\end{proof}

\begin{prop}
{Let $\rho_1,\rho_2:\pi_*\cV\to\cQ$ be two admissible quotients
and let $\cE_1=\ker\rho_1, \cE_2=\ker\rho_2$. Every isomorphism
$u:\cE_1\cong\cE_2$ can be extended to an isomorphism
$\psi:\pi_*\cV\cong\pi_*\cV$, i.e. we have a commutative diagram
\begin{equation*}
\begin{CD}
0 @>>>\cE_1@>>>\pi_*\cV@>\rho_1>>\cQ@>>> 0\\
@.       @VVuV        @VV\psi V      @VV\cong V \\
0 @>>>\cE_2@>>>\pi_*\cV@>\rho_2>>\cQ@>>> 0\\
\end{CD}
\end{equation*}
}\end{prop}

The next proposition deals with the automorphism group of
$\pi_*\cV$.

\begin{prop}
{Let $\cV$ be a locally free sheaf on $\tilde C$. Every
automorphism of $\pi_*\cV$ as an $\cO_C$-module can be induced
from an automorphism of $\cV$ as an $\cO_{\tilde C}$-module. Hence
there is a canonical isomorphism $Aut_{\cO_C}(\pi_*\cV)\cong
Aut_{\cO_{\tilde C}}(\cV)$. }\end{prop}

\begin{proof}
{Let $u:\pi_*\cV\to\pi_*\cV$ be an automorphism of $\pi_*\cV$ as
an $\cO_C$-module. It induces canonically an automorphism $\bar
u:\pi^*\pi_*\cV\to\pi^*\pi_*\cV$ as an $\cO_{\tilde C}$-module.
Let $\cT'\subset \pi^*\pi_*\cV$ be the torsion part. Then $\bar
u(\cT')=\cT'$, and it induces an automorphism $u^{\sharp}:
(\pi^*\pi_*\cV)^{\sharp}\to(\pi^*\pi_*\cV)^{\sharp}$. Since $\cV$
is locally free, there is a canonical isomorphism
$(\pi^*\pi_*\cV)^{\sharp}\cong\cV$. We obtain an automorphism
$\tilde u:\cV\to\cV$ as an $\cO_{\tilde C}$-module. Since every
step is functorial, it establishes an isomorphism
$Aut_{\cO_C}(\pi_*\cV)\cong Aut_{\cO_{\tilde C}}(\cV)$.
}\end{proof}

Next we assume $C$ is a rational nodal curve with $n$ nodes. We
describe a method to calculate $e(\M_C(r,n))$.

Let $\cE$ be a stable sheaf in $\M_C(r,n)$, and let
$\cV=(\pi^*\cE)^\sharp$ be the torsion free part of $\pi^*\cE$.
Then $\cV$ be a locally free sheaf of rank $r$ on $\tilde C$.
Since $\tilde C\cong\mathbf P^1$, by Grothendieck's Lemma,
$\cV\cong\cO(l_1)\oplus\cO(l_2)\oplus\cdots\oplus\cO(l_r)$ for
some integers $l_1\le l_2\le\cdots\le l_r$. There is a
decomposition of $\M_C(r,n)$,
\[\M_C(r, n)=\bigsqcup \M^{l_1,\cdots,l_r}_{a_1, a_2,\cdots, a_n},\]
such that $[\cE]\in \M^{l_1,\cdots,l_r}_{a_1, a_2, \cdots, a_n}$
if and only if
\[(\pi^*\cE)^\sharp\cong\cO(l_1)\oplus\cO(l_2)\oplus\cdots\oplus\cO(l_r)\]  and
\[\cE_{x_i}\cong\cO_{x_i}^{\oplus a_i}\oplus m_{x_i}^{\oplus
(r-a_i)}.\]

Let $\cE_1$, $\cE_2$ be the kernels of two admissible quotients
$\rho_1,\rho_2:\pi_*\cV\to\cQ$ respectively. The automorphism
group of $\cQ$ is a direct sum of automorphism groups of
$\cQ_{x_i}$. Let $G_i=Aut(\cQ_i)$. Then $G_i\cong GL(a_i, \mathbf
C)$. There is an $Aut(\cV)\times\prod G_i$ action on
$\Hom(\pi_*\cV, \cQ)$, $\rho\map g\circ\rho\circ u$, where
$\rho\in \Hom(\pi_*\cV, \cQ)$, $u\in Aut(\cV)$ and $g\in\prod
G_i$. Proposition 4.3 says that $\cE_1\cong\cE_2$ if and only if
$\rho_1$ and $\rho_2$ lie in the same orbit of $\Hom(\pi_*\cV,
\cQ)$ under this group action.

Next we work out a matrix form of these results under suitable
bases.

Let $V_i$ and $W_i$ be the fibres of $\cV$ at $p_i$ and $q_i$
respectively. Then $(\pi_*\cV)\otimes \mathbf C_{x_i}\cong
V_i\oplus W_i$. Since $\cQ_{x_i}=\mathbf C^{\oplus a_i}$, every
homomorphism $\rho:\pi_*\cV\to\cQ$ gives an element in the vector
space
\[ U=\oplus_{i=1}^n(\Hom(V_i, \mathbf C^{\oplus a_i})\oplus \Hom(W_i,
\mathbf C^{\oplus a_i})). \]

Fix an isomorphism
$\cV\cong\cO(l_1)\oplus\cO(l_2)\oplus\cdots\oplus\cO(l_r)$ once
and for all. For any summand $\cO(l_i)$, there is an isomorphism
of stalks $\cO(l_i)_x\cong\cO_x$ by the locally freeness of
$\cO(l_i)$. Those isomorphisms at $p_i$ and $q_i$ give rise to
bases $e_i^k\in V_i$ and $f_i^k\in W_i$. Fix all these choices
once and for all. Now an element $\rho\in U$ corresponds to a set
of $a_i\times r$ matrixes
\[\{A_i, B_i\}_{i=1, 2, \cdots, n}.\]
Let $\rho'_i\in Hom(V_i, \mathbf C^{\oplus a_i}), \rho''_i\in
Hom(W_i, \mathbf C^{\oplus a_i})$ and $v_i=(v_i^1, v_i^2, \cdots,
v_i^r)^t\in V_i$, $w_i=(w_i^1, w_i^2, \cdots, w_i^r)^t\in W_i$.
Then
\[\rho'_i(v_i)=A_i\left(\begin{matrix} v_i^1 \\ v_i^2 \\\cdots \\
v_i^r
\end{matrix}\right),\rho''_i(w_i)=B_i\left(\begin{matrix} w_i^1 \\ w_i^2 \\\cdots \\
w_i^r
\end{matrix}\right).\]

\begin{cor}
{A quotient $\{A_i, B_i\}$ is admissible if and only if the ranks
of $A_i$ and $B_i$ are both equal to $a_i$ for all $i$. In
particular, for an admissible quotient $\{A_i, B_i\}$, one has
$a_i\le r$. }\end{cor}

\begin{proof}
{Follows from proposition 4.2. }\end{proof}

Now we consider the $Aut(\cV)\times\prod G_i$ action on
$\Hom(\pi_*\cV, \cQ)$.

Evaluated at a closed point $x\in\tilde C$, every automorphism
$u\in Aut(\cV)$ gives rise to an automorphism in $Aut(V_x)$, where
$V_x$ is the fibre of $\cV$ at $x$. Therefore, every $u\in
Aut(\cV)$ gives rise to an element
\[\prod_i u(p_i)\times\prod_i
u(q_i)\in \prod_i Aut(V_i)\times\prod_i Aut(W_i).\] Let
$G'\subset\prod_i Aut(V_i)\times\prod_i Aut(W_i)$ be the subgroup
of elements derived in this way. Since $G_i\cong GL(a_i, \mathbf
C)$, there is an $G'\times\prod GL(a_i, \mathbf C)$ action on the
vector space $U$ of all quotients $\{A_i, B_i\}_{i=1, 2, \cdots,
n}$, which is given by
\[\{A_i, B_i\}\map\{g_iA_iu(p_i), g_iB_iu(q_i)\},\]
where $\prod u(p_i)\times\prod u(q_i)\in G'$ and $g_i\in GL(a_i,
\mathbf C)$. Two quotients $\{A_i, B_i\}$ and $\{A'_i, B'_i\}$ are
equivalent if they lie in one and the same orbit under this group
action.

An application of this formulation is to determine whether
$\cE\otimes\cL\cong\cE$ for $\cL\in \Pic^0C$.

Let $\rho:\pi_*\cV\to\cQ$ be an admissible quotient with
corresponding matrixes $\{A_i, B_i\}$ and let $\cE=\ker\rho$. Let
$\cL\in \Pic^0 C$ be given by the matrixes $\{1, t_i\}$, where
$t_i\in\mathbf C^*$. The exact sequence
\[0\map\cE\map\pi_*\cV\stackrel{\rho}{\map}\cQ\map 0\] induces
an exact sequence
\[0\map\cE\otimes\cL\map\pi_*\cV\otimes\cL\stackrel{\rho\otimes 1}{\map}\cQ\otimes\cL\map 0.\]
Note that $\pi_*\cV\otimes\cL\cong\pi_*\cV$ and the quotient
$\rho\otimes 1:\pi_*\cV\otimes\cL\to\cQ\otimes\cL$ is also
admissible. Fix an isomorphism $\cQ\otimes\cL\cong\cQ$ and choose
corresponding bases, $\rho\otimes 1$ is given by the matrixes
$\{A_i, t_iB_i\}$.

We are now ready to calculate the Euler number $e(\M_C(r, n))$.
For the purpose of this paper, we only consider the case $r=2$ and
$n=1$.

\begin{prop}
{Let $\M^{l_1,l_2}_{a_1, a_2, \cdots, a_n}$ be a stratum in
$\M_C(2, 1)$ such that $\sum a_i\ge 2$. Then $e(\M^{l_1,l_2}_{a_1,
a_2, \cdots, a_n})=0$. }\end{prop}

\begin{proof}
{For simplicity, we consider only the stratum $\M^{l_1,l_2}_{1, 1,
0, \cdots, 0}$ as illustration. Let $[\cE]\in \M^{l_1,l_2}_{1, 1,
0, \cdots, 0}$ be the kernel of an admissible quotient $\{A_i,
B_i\}$. We can choose a suitable base such that $A_1=A_2=(1, 0),
B_1=B_2=(0, 1)$. For an odd prime $p$, let $\cL$ be given by
$t_1=1, t_2=\zeta$, where $\zeta$ is a $p$-th primitive root of
unity. Then $\cL^{\otimes p}=\cO_C$ and $\cE\otimes\cL$
corresponds to the quotient $\{A_i, t_iB_i\}$. It is direct to
verify that $\{A_i, t_iB_i\}$ and $\{A_i, B_i\}$ are not
equivalent, hence $\cE\otimes\cL$ and $\cE$ are not isomorphic. So
we get a free $\mathbf Z/(p)$ action on $\M^{l_1,l_2}_{1, 1, 0,
\cdots, 0}$. Because $p$ can be chose arbitrarily large,
$e(\M^{l_1,l_2}_{1, 1, 0, \cdots, 0})=0$. }\end{proof}

Since $\M^{l_1,l_2}_{0, 0, \cdots, 0}$ is empty, by this
proposition, the contribution to the Euler number $e(\M_C(2, 1))$
comes from strata $\M^{l_1,l_2}_{a_1, a_2, \cdots, a_n}$ with
$\sum a_i=1$. Because $\chi(\cE)=1$ and $\cE$ fits into an exact
sequence
\[ 0\map\cE\map\pi_*(\cO(l_1)\oplus\cO(l_2))\map\mathbf C_{x_i}\map 0,\]
the stability of $\cE$ forces $l_1=l_2=0$. Every $\M^{0,0}_{0,
\cdots, 1,\cdots, 0}$ is a set of a single point. Therefore,
\begin{prop}
{The Euler number $e(\M_C(2, 1))$ is equal to $n$, which is the
number of nodes on $C$. }\end{prop}

Let $D=2C_0$ be a divisor in the set $\cW_3^0$. Then $\M_D^1$ is
isomorphic to $\M_{C_0}(2,1)$. The number of nodes on $C_0$ is
equal to the arithmetic genus $g=\frac{1}{2}H^2+1$ of $C_0$.
Therefore,
\begin{prop}
{$e(\frak M_D^1)=g$. }\end{prop}

\section {Calculation of $e(\Phi^{-1}(\cW_3^0))$, Part II}

This is the second part of the calculation of
$e(\Phi^{-1}(\cW_3^0))$. As we mentioned in the previous section,
$\M_D$ is a disjoint union of $\M_D^1$ and $\M_D^2$ for $D\in
\cW_3^0$. We have calculated $e(\M_D^1)$. In this section, we will
show that $e(\M_D^2)=0$.

Let $C_0\subset S$ be a nodal curve, and let $C$ be the associated
nonreduced curve to the divisor $2C_0$. Let $p$ be a node on
$C_0$, and $\pi_0: \hat {C_0} \to C_0$ be the partial
normalization of $C_0$ at $p$. Now we construct a curve $\hat C$,
which is an infinitesimal extension of $\hat{C_0}$, and a finite
morphism $\pi:\hat C\to C$ called a partial normalization of $C$.

We pick a small neighborhood $U$ around $p$ on the surface $S$,
such that $C$ is defined by $x^2y^2=0$ in $U$. Let $\mathbf C\{x,
y\}$ be the ring of holomorphic functions on $U$. Then
$\cO_C(U\cap C)=\mathbf C\{x, y\}/(x^2y^2)$. The injective
homomorphism
\[\psi:\mathbf C\{x, y\}/(x^2y^2)\to \mathbf C\{x, u\}/(u^2)\oplus
\mathbf C\{v, y\}/(v^2)\] is a local isomorphism except at $p$.
Remove the point $p$ on $C$ and glue the pieces defined by the
ringed space $\mathbf C\{x, u\}/(u^2)\oplus \mathbf C\{v,
y\}/(v^2)$ along $\psi$, we get a curve $\hat C$, and a finite map
$\pi:\hat C\to C$. There is a canonical exact sequence \[ 0\map
\cO_C\map\pi_*\cO_{\hat C}\map\cA\map 0,
\] where $\cA\cong \mathbb C[x,y]/(x^2,y^2)$ is a skyscraper sheaf
supported at $p$. Moreover, there is a commutative diagram
\begin{equation*}
\begin{CD}
0 @>>>\cO_C @>>>\pi_*\cO_{\hat C} @>>>\cA @>>> 0   \\
@.   @VV V          @VV V       @V VV     \\
0 @>>>\cO_{C_0} @>>>\pi_*\cO_{\hat C_0} @>>>\mathbf C_p @>>> 0  \\
\end{CD}
\end{equation*}

Let $\cI$ and $\hat\cI$ be the nilpotent ideal sheaves of $\cO_C$
and $\cO_{\hat C}$ respectively. Then $\chi(\hat\cI)=\chi(\cI)+3$.
It implies that $\deg\hat\cI=\deg\cI+2$. Let $C$ be a rational
nodal curve on a $K3$ surface and let $\tilde C\to C$ be the
normalization of $C$. Then $\deg\cI=2-2g=H^2$. Because the number
of nodes on $C$ is equal to $g$, $\deg\tilde\cI=\deg\cI+2g=2$.

\begin{prop}
{Let $C$ be a nonreduced curve with nilpotent ideal sheaf $\cI$.
Suppose $\cI$ is invertible as a sheaf of $\cO_{C_0}$-modules and
$\deg\cI>0$. Let $\cE$ be a pure sheaf of $\cO_C$-modules such
that $\cE_{\eta}\cong\cO_{\eta}$ at the generic point $\eta$ of
$C$. Then $\cE$ is not stable.

}\end{prop}

\begin{proof}
{Let $\cE$ be such a sheaf and let $\cE^{\sharp}_0$ be the torsion
free part of $\cE_0=\cE\otimes\cO_{C_0}$, considered as a sheaf of
$\cO_{C_0}$-modules. There is a canonical homomorphism
$\cE\to\cE_0^{\sharp}$. Every quotient $\cE\to\cF$ with $\cF$ a
torsion free $\cO_{C_0}$-module is equivalent to $\cE\to
\cE_0^{\sharp}$. Therefore, for the stability of $\cE$, it is
enough to check the quotient $\cE\to \cE_0^{\sharp}$.

We start with the exact sequence
\[0\map\cI\map\cO_C\map\cO_{C_0}\map 0.\]
Tensoring with $\cE$, we obtain
\[\cE_0\otimes\cI\map\cE\map\cE_0\map 0.\]
Let $\cT'$ be the torsion part of $\cE_0\otimes\cI$, and let
$(\cE_0\otimes\cI)^{\sharp}=(\cE_0\otimes\cI)/\cT'$. There is an
exact sequence
\[0\map(\cE_0\otimes\cI)^{\sharp}\map\cE\map \cE_0\map 0.\]
On the other hand, we have
\[0\map\cK\map\cE\map \cE_0^{\sharp}\map 0\]
Because $\cI$ is an invertible sheaf of $\cO_{C_0}$-modules, the
torsion part of $\cE_0$ is isomorphic to $\cT'$.
$\chi(\cK)=\chi(\cE)-\chi(\cE_0^{\sharp})=\chi(\cE_0\otimes\cI)$.
Because
$\chi(\cE_0\otimes\cI)=\chi(\cE_0)+\deg\cI>\chi(\cE_0)>\chi(\cE_0^{\sharp})$,
$\cE$ is not stable. }\end{proof}

Let $\pi:\hat C\to C$ be the partial normalization of $C$ at $p$.
Let $\cF$ be a sheaf of $\cO_C$-modules which is pure of dimension
$1$. Then there is a canonical homomorphism
$\cF\to\pi_*(\pi^*\cF)$. Let $T_0\subset \pi^*\cF$ be the maximal
subsheaf of dimension $0$, we get a sheaf $(\pi^*
\cF)^{\sharp}=\pi^*\cF/\cT_0$ which is pure of dimension 1, and
there is a canonical injective homomorphism $\cF\to\pi_*(\pi^*
\cF)^{\sharp}$. The cokernel $\cT$ is a skyscraper sheaf supported
at $p$. Note that if $\cF$ satisfies $\cI\cF\ne 0$, so does
$(\pi^* \cF)^{\sharp}$ as a sheaf of $\cO_{\hat C}$-modules.

The notion of admissible quotients can be defined in the same way
as in section 4, and propositions 4.2-4.4 are also true in this
case.

Let $\rho:\pi_*\cE\to\cQ$ be an admissible quotient, and let $p$
be a node on $C$ with $\pi^{-1}(p)=\{q_1,q_2\}$. Then
$(\pi_*\cE)_p=\cE_{q_1}\oplus \cE_{q_2}$. We let $\iota_i:
\cE_{q_i}\to(\pi_*\cE)_p$ and $p_i: (\pi_*\cE)_p\to\cE_{q_i}$ be
the natural injections and projections respectively. Define
$\rho^i$ as compositions
\[\rho^i:
\cE_{q_i}\stackrel{\iota_i}{\to}(\pi_*\cE)_p\stackrel{\rho}{\to}
\cQ_p.\] Because $\rho$ is admissible, $\rho^i$ are both
surjective homomorphisms. Clearly $\rho=\rho^1p_1+\rho^2p_2$. For
$t\in \mathbf C^*$, we define $\rho_t=\rho^1p_1+t\rho^2p_2$. It
gives rise to a surjective homomorphism $\rho_t: \pi_*\cE\to \cQ$
which is also admissible.

Let $\pi:\hat C\to C$ be the partial normalization of $C$ at $p$.
Apply the above construction to the canonical exact sequence
\[0\map \cO_C\map \pi_*\cO_{\hat C}\stackrel{\rho}{\map}\cA\map
0,\] and let $\cL_t=\ker\rho_t$. Then $\cL_t$ is invertible for
every $t\in \mathbf C^*$. The set of these invertible sheaves form
a subgroup $G_p\subset \Pic^0 C$. Clearly $G_p\cong \mathbf C^*$.

For any admissible quotient $\rho: \pi_*\cE\to\cQ$, let $\cK_t$ be
the kernel of $\rho_t$.

\begin{lem}
There is a canonical isomorphism between $\cL_t\otimes\cK_s$ and
$\cK_{st}$.
\end{lem}

Now we give a decomposition on $\M_D^2$ for $D\in \cW_3^0$. Let
$\pi:\tilde C\to C$ be the normalization of $C$. Let
$\M_{\tilde\cF,\cT}$ be the subset consists of stable sheaves
$\cF$ such that $(\pi^*\cF)^\sharp\cong\tilde\cF$ and
$\pi_*\tilde\cF/\cF\cong\cT$. We get a decomposition
$\M_D^2=\sqcup \M_{\tilde\cF,\cT}$. In fact, for every nonempty
stratum $\M_{\tilde\cF,\cT}$, $\cT$ is nonzero. Because if
$\M_{\tilde\cF,0}$ is nonempty, a sheaf $\cF$ in
$\M_{\tilde\cF,0}$ is the direct image of a sheaf on $\tilde C$,
i.e. $\cF\cong\tilde\pi_*\cE$ for a sheaf $\cE$ of $\cO_{\tilde
C}$-modules. By proposition 5.1, $\cE$ is not stable, which
violates the stability of $\cF$.

Let $\M_{\tilde\cF,\cT}$ be a stratum, and let $p\in C$ be a node
such that $\cT_p\ne 0$. There is a subgroup $G_p\subset \Pic^0 C$
defined as above, and a $G$-action on $\M_{\tilde\cF,\cT}$ defined
by tensorization. Next we will show that this group action is free
on $\M_{\tilde\cF,\cT}$. The following lemma is useful in the
proof.

Let $\cE$ be a pure sheaf of $\cO_C$-modules such that
$\cE_{\eta}\cong\cO_{\eta}$ at the generic point $\eta$ of $C$.
Let $\cE''$ be the torsion free part of
$\cE_0=\cE\otimes_{\cO_C}\cO_{C_0}$, and let $\cE'$ be the kernel
of the restriction homomorphism $f:\cE\to \cE''$. Since $f$ is not
an isomorphism, $\cE'$ and $\cE''$ are both rank $1$ torsion free
sheaves of $\cO_{C_0}$-modules whose automorphism groups are
$\mathbf C^*$.
\begin{lem}
{Let $\psi:\cE\to\cE$ be an automorphism and let $c: \cE'\to\cE'$
and $d: \cE''\to\cE''$ be the induced automorphisms. Then they fit
into the commutative diagram
\begin{equation*}
\begin{CD}
0 @>>> \cE' @>>> \cE @>>> \cE'' @>>> 0\\
@.       @VVcV          @VV\psi V            @VVdV \\
0 @>>> \cE' @>>> \cE @>>> \cE'' @>>> 0\\
\end{CD}
\end{equation*}
and $c=d$.
}\end{lem}

\begin{proof}
{Consider $\psi'=\psi-c\cdot id:\cE\to\cE$. Clearly
$\psi'(\cE')=0$. It induces a homomorphism $u:\cE''\to\cE$.
Composed with $\cE\to\cE''$, we get $h:\cE''\to\cE''$. Since
$\cE''$ is torsion free and has rank $1$ as an $\cO_{C_0}$-module,
$h$ is a multiplication by $(d-c)$. If $h\ne 0$, then after
scaling $u$ by $\frac{1}{h}$, $u$ splits the exact sequence and
hence $\cE\cong\cE'\oplus\cE''$, which contradicts to
$\cE_{\eta}\cong\cO_{\eta}$. Therefore $h=0$, i.e. $c=d$.
}\end{proof}

\begin{prop}
{Let $\M_{\tilde\cF,\cT}$ be a stratum such that $\cT_p\ne 0$ for
a node $p\in C$. Then the associated $G_p$-action on
$\M_{\tilde\cF,\cT}$ is free. Therefore, by the decomposition of
$\M_D^2$, $e(\M_D^2)=0$. }\end{prop}

\begin{proof}
{Let $\pi:\hat C\to C$ be the partial normalization of $C$ at $p$.
Let $\cF$ be a stable sheaf in $\M_{\tilde\cF,\cT}$. Then $\cF$
fits into the exact sequence
\[0\map\cF\map\pi_*(\pi^*\cF)^{\sharp}\map\cT\map 0.\]
Let $\cE=(\pi^*\cF)^{\sharp}$. Then $\rho:\pi_*\cE\map\cT$ is
clearly an admissible quotient. We let $\cK_t$ be the kernel of
$\rho_t$. Then $\cF=\cK_1$ and $\cF\otimes\cL_t=\cK_t$. Suppose
$\cF\otimes\cL_t\cong \cF$ for some $\cL_t\in G_p$, there is a
commutative diagram
\begin{equation*}
\begin{CD}
0 @>>> \cK_1 @>>>\pi_*\cE@>\rho_1>>\cT@>>> 0\\
@. @VV\cong V     @VV\psi V      @VVhV \\
0 @>>> \cK_t @>>>\pi_*\cE@>\rho_t>>\cT@>>> 0\\
\end{CD}
\end{equation*}

It induces the following diagram on the stalks at the node $p$,
\begin{equation*}
\begin{CD}
\cE_{q_1}\oplus\cE_{q_2}@>\rho_1>>\cT_p@>>> 0\\
          @VV\psi V      @VVhV \\
\cE_{q_1}\oplus\cE_{q_2}@>\rho_t>>\cT_p@>>> 0\\
\end{CD}
\end{equation*}
where $\{q_1, q_2\}=\pi^{-1}(p)$.

Recall that there is a canonical exact sequence for $\cE$,
\[0\map\cE'\map\cE\map\cE''\map 0,\] where $\cE'$ and
$\cE''$ are nonzero torsion free sheaves of $\cO_{C_0}$-modules.
Let $\cT^i\subset\cT_p$ be the images of $\cE'_{q_i}$ under the
surjective homomorphisms $\rho^i:\cE_{q_i}\to\cT_p$. Then we have
the following three cases.

Case 1. $\cT^1\ne 0$ and $\cT^2\ne 0$.

Since $\psi:\pi_*\cE\to \pi_*\cE$ is induced from an automorphism
of $\cE$, $\psi(\cE_{q_i})=\cE_{q_i}$. Consider the restriction of
the diagram to $\cE'_{q_i}$ respectively, by Lemma 5.3, there are
commutative diagrams
\begin{equation*}
\begin{CD}
\cE'_{q_1}@>\rho^1>>\cT^1 @. \qquad\qquad \cE'_{q_2}@>\rho^2>>\cT^2\\
          @VVc V      @VVh_1V  \qquad\qquad @VVc V      @VVh_2V \\
\cE'_{q_1}@>\rho^1>>\cT^1 @. ,\qquad\qquad \cE'_{q_2}@>t\rho^2>>\cT^2.\\
\end{CD}
\end{equation*}
If $\cT^1\cap\cT^2\ne\{0\}$, let $0\ne x\in \cT^1\cap\cT^2$. Then
from the left diagram, $h(x)=h_1(x)=cx$, and from the right
diagram, $h(x)=h_2(x)=ctx$. It implies that $t=1$ and therefore
the group action is free. Next we assume $\cT^1\cap\cT^2=\{0\}$.
Let $x\in \cT^1$ be a nonzero element. Then the image $\bar x$ of
$x$ in $\cT_p/\cT^2$ is nonzero. From the commutative diagram
\begin{equation*}
\begin{CD}
\cE''_{q_2}@>\rho^2>>\cT_p/\cT^2\\
          @VVc V      @VVh_3V \\
\cE''_{q_2}@>t\rho^2>>\cT_p/\cT^2\\
\end{CD}
\end{equation*}
we have $h_3(\bar x)=ct\bar x$. Because $h(x)=h_1(x)=cx$,
$h_3(\bar x)=c\bar x$. It implies that $t=1$.

Case 2. $\cT^1\ne 0$ and $\cT^2=0$(Or equivalently $\cT^1=0$ and
$\cT^2\ne 0$).

Since $\cT^2=0$ and $\cE_{q_2}\stackrel{\rho^2}{\to}\cT_p$ is
surjective, there is a surjective morphism $\rho^2:
\cE''_{q_2}\to\cT_p$ and a commutative diagram
\begin{equation*}
\begin{CD}
\cE''_{q_2}@>\rho^2>>\cT_p\\
          @VVc V      @VVhV \\
\cE''_{q_2}@>t\rho^2>>\cT_p\\
\end{CD}
\end{equation*}
Let $0\ne x\in \cT^1$. We already have $h(x)=h_1(x)=cx$. However,
from the above diagram, $h(x)=ctx$. Hence $t=1$.

Case 3. $\cT^1=0$ and $\cT^2=0$.

Since $\cT^1=0$ and $\cT^2=0$, we get commutative diagrams
\begin{equation*}
\begin{CD}
\cE''_{q_1}@>\rho^1>>\cT_p @. \qquad\qquad \cE''_{q_2}@>\rho^2>>\cT_p\\
          @VVc V      @VVh_1V  \qquad\qquad @VVc V      @VVh_2V \\
\cE''_{q_1}@>\rho^1>>\cT_p @. ,\qquad\qquad \cE''_{q_2}@>t\rho^2>>\cT_p \\
\end{CD}
\end{equation*}
Apply the same argument as in case 1, we get $t=1$.}\end{proof}

Because the cardinality of the finite set $\cW_3^0$ is $G_g$,
combine proposition 4.8 and 5.4, we conclude

\begin{prop}
{$e(\Phi^{-1}(\cW_3^0))=gG_g$.

}\end{prop}

\end{document}